\documentstyle[12pt]{article}
\begin{document}
\title {On the Orthogonal Stability of the Pexiderized Quadratic Equation }
\author{{\bf Mohammad Sal Moslehian} \\ Dept. of Math., Ferdowsi Univ.\\ P. O. Box 1159, Mashhad 91775\\ Iran\\ E-mail: msalm@math.um.ac.ir\\Home: http://www.um.ac.ir/$\sim$moslehian/}
\date{}
\maketitle
\begin{abstract}
The Hyers--Ulam stability of the conditional quadratic functional equation of Pexider type $f(x+y)+f(x-y)=2g(x)+2h(y), x\perp y$ is established where $\perp$ is a symmetric orthogonality in the sense of R\" atz and $f$ is odd.
\end{abstract}
\newpage
\section{Introduction.}

Let us denote the sets of real and nonnegative real numbers by ${\bf R}$ and ${\bf R_+}$, respectively. 

Suppose that $X$ is a real vector space with $\dim X\geq 2$ and $\perp$ is a binary relation on $X$ with the following properties:\\
(O1) totality of $\perp$ for zero: $x\perp 0, 0\perp x$ for all $x\in X$;\\
(O2) independence: if $x,y\in X-\{0\}, x\perp y$, then $x,y$ are linearly independent;\\
(O3) homogeneity: if $x,y\in X, x\perp y$, then $\alpha x\perp \beta y$ for all $\alpha, \beta\in {\bf R}$;\\
(O4) the Thalesian property: Let $P$ be a $2$-dimensional subspace of $X$. If $x\in P$ and $\lambda\in{\bf R_+}$, then there exists $y_0\in P$ such that $x\perp y_0$ and $x+y_0\perp \lambda x-y_0$.

The pair $(X,\perp)$ is called an orthogonality space. By an orthogonality normed space we mean an orthogonality space equipped with a norm. 

Some examples of special interest are (i) The trivial orthogonality on a vector space $X$ defined by (O1), and for non-zero elements $x,y\in X$, $x\perp y$ if and only if $x,y$ are linearly independent.\\
(ii) The ordinary orthogonality on an inner product space $(X, \langle.,.\rangle)$ given by $x\perp y$ if and only if $\langle x,y\rangle=0$.\\
(iii) The Birkhoff-James orthogonality on a normed space $(X,\|.\|)$ defined by $x\perp y$ if and only if $\|x+\lambda y\|\geq \|x\|$ for all $\lambda\in {\bf R}$.

The relation $\perp$ is called symmetric if $x\perp y$ implies that $y\perp x$ for all $x,y\in X$. Clearly examples (i) and (ii) are symmetric but example (iii) is not. It is remarkable to note, however, that a real normed space of dimension greater than or equal to $3$ is an inner product space if and only if the Birkhoff-James orthogonality is symmetric.

Let $X$ be a vector space (an orthogonality space) and $(Y,+)$ be an abelian group. A mapping $f:X\to Y$ is called (orthogonally) additive if it satisfies the so-called (orthogonal) additive functional equation $f(x+y)=f(y)+f(x)$ for all $x,y\in X$ (with $x\perp y$). A mapping $f:X\to Y$ is said to be (orthogonally) quadratic if it satisfies the so-called (orthogonally) Jordan-von Neumann quadratic function equation $f(x+y)+f(x-y)=2f(x)+2f(y)$ for all $x,y\in X$ (with $x\perp y$). For example, a function $f:X\to Y$ between real vector spaces is quadratic if and only if there exists a (unique) symmetric bi-additive mapping $B:X\times X\to Y$ such that $f(x)=B(x,x)$ for all $x\in X$. In fact, $B(x,y)=\frac{1}{4}(f(x+y)-f(x-y))$, cf.  $\cite{A-D}$.

In the recent decades, stability of functional equations have been investigated by many mathematicians. They have so many applications in Information Theory, Economic Theory and Social and Behaviour Sciences; cf. $\cite{ACZ}$.

The first author treating the stability of the quadratic equation was F. Skof $\cite{SKO}$ by proving that if $f$ is a maping from a normed space $X$ into a Banach space $Y$ satisfying $\|f(x+y)+f(x-y)-2f(x)-2f(y)\|\leq \epsilon$ for some $\epsilon>0$, then there is a unique quadratic function $g:X\to Y$ such that $\|f(x)-g(x)\|\leq\frac{\epsilon}{2}$. P. W. Cholewa $\cite{CHO}$ extended Skof's theorem by replacing $X$ by an abelian group G. Skof's result was later generalized by S. Czerwik $\cite{CZE1}$ in the spirit of Hyers--Ulam--Rassias. K.W. Jun and Y. H. Lee $\cite{J-L}$ proved the stability of quadratic equation of Pexider type. The stability problem of the quadratic equation has been extensively investigated by some mathematicians $\cite{RAS}, \cite{CZE2}, \cite{CZE3}$.

The orthogonal quadratic equation $$f(x+y)+f(x-y)=2f(x)+2f(y),~ x\perp y$$
was first investigated by F. Vajzovi\' c $\cite{VAJ}$ when $X$ is a Hilbert space, $Y$ is the scalar field, $f$ is continuous and $\perp$ means the Hilbert space orthogonality. Later H. Drljevi\' c $\cite{DRL}$, M. Fochi $\cite{FOC}$ and G. Szab\' o $\cite{SZA}$ generalized this result. 

One of the significant conditional equations is the so-called orthogonally quadratic functional equation of Pexider type $f(x+y)+f(x-y)=2g(x)-2h(y),~~ x\perp y$. Our main aim is to consider the stability of this equation in the spirit of Hyers--Ulam under certain conditions.

\section {Main Result.}

Applying some ideas from $\cite{G-S}$ and $\cite{J-S}$, we deal with the conditional stability problem for 
$$f(x+y)+f(x-y)=2g(x)+2h(y)~~~~~x\perp y$$
where $f$ is odd and $\perp$ is symmetric. We also use a sequence of Hyers' type $\cite{HYE}$ which is a useful tool in the theory of  stability of equations. (See also $\cite{MOS}$)

Throughout this section, $(X,\perp)$ denotes an orthogonality normed space and $(Y, \|.\|)$ is a real Banach space.

{\bf 2.1. Lemma.} If $A:X\to Y$ fulfills $A(x+y)+A(x-y)=2A(x)$ for all $x,y\in X$ with $x\perp y$ and $\perp$ is symmetric, then $A(x)-A(0)$ is orthogonally additive.

{\bf Proof.} Assume that $A(x+y)+A(x-y)=2A(x)$ for all $x,y\in X$ with $x\perp y$. Putting $x=0$ , we get $-A(y)=A(-y)-2A(0), y\in X$. Let $x\perp y$. Then $y\perp x$ and so $A(y-x)=-A(y+x)+2A(y)$. Hence $A(x+y)=-A(x-y)+2A(x)=(A(y-x)-2A(0))+2A(x)=(-A(y+x)+2A(y))-2A(0)+2A(x)$.\\
Thus $A(x+y)-A(0)=(A(x)-A(0))+(A(y)-A(0))$. So that $A(x)-A(0)$ is orthogonally additive.$\Box$

{\bf 2.2. Remark.} Thanking J. R\" atz, there exists an odd mapping $A$ from an orthogonality space into a uniquely 2-divisible group $(Y,+)$ (i.e. an abelian group in which the map $\varphi:Y\to Y, \varphi(x)=2x$ is bijective) satisfying $A(x+y)+A(x-y)=2A(x), x\perp y$ such that $A(0)\neq 0$. Consider $Y={\bf Z}_2=\{0,1\}$ and $A(x)=1, x\in X$.

{\bf 2.3. Theorem.} {\it Suppose $\perp$ is symmetric on $X$ and $f, g, h:X\to Y$ are mappings fulfilling 
\begin{equation}
\|f(x+y)+f(x-y)-2g(x)-2h(y)\|\leq\epsilon
\end{equation}
for some $\epsilon$ and for all $x,y\in X$ with $x\perp y$. Assume that $f$ is odd. Then there exist exactly an additive mapping $T:X\to Y$ and exactly a quadratic mapping $Q:X\to Y$ such that
\begin{eqnarray*}
\|f(x)-T(x)-Q(x)\|\leq 6\epsilon
\end{eqnarray*}
\begin{eqnarray*}
\|g(x)-T(x)-Q(x)\|\leq \frac{13}{2}\epsilon
\end{eqnarray*}
for all $x\in X$.}

{\bf Proof.} Put $x=0$ in $(1)$. We can do this because of (O1). Then
\begin{equation}
\|h(y)\|\leq\frac{1}{2}\epsilon
\end{equation}
for all $y\in X$. Similarly, by putting $y=0$ in $(1)$ we get
\begin{equation}
\|f(x)-g(x)\|\leq\frac{1}{2}\epsilon
\end{equation}
for all $x\in X$. Hence
\begin{eqnarray}
\|f(x+y)+f(x-y)-2f(x)\|&\leq& \|f(x+y)+f(x-y)-2g(x)-2h(y)\|\nonumber\\
&+& 2\|f(x)-g(x)\|+2\|h(y)\| \nonumber\\
&\leq&\epsilon+2.\frac{1}{2}\epsilon+2.\frac{1}{2}\epsilon=3\epsilon 
\end{eqnarray}
for all $x,y\in X$ with $x\perp y$.

Fix $x\in X$. By (O4), there exists $y_0\in X$ such that $x\perp y_0$ and $x+y_0\perp x-y_0$. Since $\perp$ is symmetric $x-y_0\perp x+y_0$, too. Using inequality $(4)$ and the oddness of $f$ we get
\begin{eqnarray*}
\|f(x+y_0)+f(x-y_0)-2f(x)\|\leq 3\epsilon\\
\|f(2x)+f(2y_0)-2f(x+y_0)\|\leq 3\epsilon\\
\|f(2x)-f(2y_0)-2f(x-y_0)\|\leq 3\epsilon.
\end{eqnarray*} 
So that
\begin{eqnarray*}
&&\|f(2x)-2f(x)\|\leq \|f(x+y_0)+f(x-y_0)-2f(x)\|\\
&+&\frac{1}{2}\|f(2x)+f(2y_0)-2f(x+y_0)\|+\frac{1}{2}\|f(2x)-f(2y_0)-2f(x-y_0)\|\\
&\leq&6\epsilon.
\end{eqnarray*}
It is not hard to see that 
\begin{eqnarray*}
\|2^{-n}f(2^nx)-f(x)\|\leq 6\epsilon\displaystyle{\sum_{k=1}^{n}}(\frac{1}{2})^k
\end{eqnarray*}
for all $n$, and 
\begin{eqnarray*}
\|2^{-n}f(2^nx)-2^{-m}f(2^mx)\|\leq 6\epsilon\displaystyle{\sum_{k=m+1}^{n}}(\frac{1}{2})^k.
\end{eqnarray*}
for all $m<n$. Thus $\{2^{-n}f(2^nx)\}$ is a Cauchy sequence in the Banach space $Y$. Hence $\displaystyle{\lim_{n\to\infty}}2^{-n}f(2^nx)$ exists and well defines the odd mapping $A(x):=\displaystyle{\lim_{n\to\infty}}2^{-n}f(2^nx)$ from $X$ into $Y$ satisfying 
\begin{equation}
\|f(x)-A(x)\|\leq 6\epsilon~~~~~x\in X.
\end{equation}
For all $x,y\in X$ with $x\perp y$, by applying inequality $(4)$ and (O3) we obtain
\begin{eqnarray*}
\|2^{-n}f(2^n(x+y))+2^{-n}f(2^n(x-y))-2^{-n+1}f(2^nx)\|\\
\leq 3.2^{-n}\epsilon.
\end{eqnarray*}
If $n\to\infty$ then we deduce that $A(x+y)+A(x-y)-2A(x)=0$ for all $x,y\in X$ with $x\perp y$. Moreover, $A(0)=\displaystyle{\lim_{n\to\infty}}2^{-n}f(2^n.0)=0$. Using Lemma 2.1 we conclude that $A$ is an ortogonally additive mapping.

By Corollary 7 of $\cite{RAT}$, $A$ therefore is of the form $T+Q$ with $T$ additive and $Q$ quadratic. If there are another quadratic mapping $Q'$ and another additive mapping $T'$ satisfying the required inequalities in our theorem and $A'=T'+Q'$ then 

\begin{eqnarray*}
\|A(x)-A'(x)\|&\leq&\|f(x)-A(x)\|+\|f(x)-A'(x)\|\leq 12\epsilon
\end{eqnarray*}
for all $x\in X$. Using the fact that additive mappings are odd and quadratic mappings are even we obtain   
\begin{eqnarray*}
\|T(x)-T'(x)\|&=&\|\frac{1}{2}[(T(x)+Q(x)-T'(x)-Q'(x))\\
&&+(T(x)-Q(x)-T'(x)+Q'(x))]\|\\
&\leq&\frac{1}{2}\|T(x)+Q(x)-T'(x)-Q'(x)\|\\
&&+\frac{1}{2}\|T(x)-Q(x)-T'(x)+Q'(x)\|\\
&\leq&\frac{1}{2}\|A(x)-A'(x)\|+\frac{1}{2}\|A(-x)-A'(-x)\|=\\
&\leq& 12\epsilon.
\end{eqnarray*}
Hence 
\begin{eqnarray*}
\|T(x)-T'(x)\|=\frac{1}{n}\|T(nx)-T'(nx)\|&\leq& \frac{12}{n}\epsilon.
\end{eqnarray*}
Tending $n$ to $\infty$ we infer that $T=T'$.
Similarly, 
\begin{eqnarray*}
\|Q(x)-Q'(x)\|&=&\|\frac{1}{2}[(T(x)+Q(x)-T'(x)-Q'(x))\\
&&-(T(x)-Q(x)-T'(x)+Q'(x))]\|\\
&\leq&\frac{1}{2}\|T(x)+Q(x)-T'(x)-Q'(x)\|\\
&&+\frac{1}{2}\|T(x)-Q(x)-T'(x)+Q'(x)\|\\
&\leq&\frac{1}{2}\|A(x)-A'(x)\|+\frac{1}{2}\|A(-x)-A'(-x)\|\\
&\leq& 12\epsilon.
\end{eqnarray*}
for all $x\in X$.
Hence 
\begin{eqnarray*}
\|Q(x)-Q'(x)\|=\frac{1}{n^2}\|Q(nx)-Q'(nx)\|\leq\frac{12}{n^2}\epsilon
\end{eqnarray*}
for all $x\in X$. Taking the limit, we conclude that $Q=Q'$.

Using $(5)$ and $(3)$ we infer that for all $x\in X$,
\begin{eqnarray*}
\|g(x)-A(x)\|\leq\|g(x)-f(x)\|+\|f(x)-A(x)\|\leq \frac{1}{2}\epsilon+6\epsilon=\frac{13}{2}\epsilon.\Box
\end{eqnarray*}

{\bf 2.4. Remark.} (i) If $g=\lambda f$ for some number $\lambda\neq 1$, then inequality $(3)$ implies that $|1-\lambda|\|f(x)\|\leq\frac{1}{2}\epsilon, x\in X$. Hence $|1-\lambda|\|2^{-n}f(2^nx)\|\leq 2^{-n-1}\epsilon, x\in X$. So $A(x)=\displaystyle{\lim_{n\to\infty}} 2^{-n}f(2^nx)=0, x\in X$. 

(ii) Similarly, if $h=\lambda f$ for some number $\lambda\neq 0$, then it follows from $(2)$ that $A(x)=0$ for all $x\in X$.

As far as the author knows, unlike orthogonally additive maps (see Corollary 7 of $\cite{RAT}$), there is no characterization for orthogonally quadratic maps. Every orthogonally quadratic mapping $q$ into a uniquely 2-divisible abelian group $(Y,+)$ is even. In fact, $0\perp 0$ and so $q(0)+q(0)=4q(0)$. Therefore $q(0)=0$. For all $y\in X$ we have $0\perp y$ and hence $q(y)+q(-y)=2q(0)+2q(y)$. Thus $q(-y)=q(y)$.

There are some characterizations of orthogonally quadratic maps in various notions of orthogonality. For example if $A$-orthogonality on a Hilbert space $H$ is defined by $\perp_A=\{(x,y): <Ax,y>=0\}$ where $A$ is a bounded self-adjoint operator on $H$ then, as shown by M. Fochi, every $A$-orthogonally quadratic functional is quadratic if $\dim A(H)\geq 3$; cf. $\cite{FOC}$ and $\cite{SZA}$.

{\bf 2.5 Problem.} Suppose $f, g, h:X\to Y$ are mappings fulfilling 
\begin{equation}
\|f(x+y)+f(x-y)-2g(x)-2h(y)\|\leq\epsilon
\end{equation}
for some $\epsilon$ and for all $x,y\in X$ with $x\perp y$. Assume that $f$ is even. Does there exist an orthogonally quadratic mapping $Q:X\to Y$, under certain conditions, such that
\begin{eqnarray*}
\|f(x)-Q(x)\|\leq \alpha\epsilon
\end{eqnarray*}
\begin{eqnarray*}
\|g(x)-Q(x)\|\leq \beta\epsilon
\end{eqnarray*}
\begin{eqnarray*}
\|h(x)-Q(x)\|\leq \gamma\epsilon
\end{eqnarray*}
for some scalars $\alpha, \beta, \gamma$ and for all $x\in X$.

{\bf Acknowledgment.} The author would like to sincerely thank Professors J. R\" atz, R. Ger and J. Sikorska for their useful comments on orthogonally quadratic mappings.


\begin{thebibliography}{99}
\bibitem{ACZ} J. Acz\' el, A short course on functional equations, D. Reidel Publ. Co., Dordrecht, 1987.
\bibitem{A-D} J. Acz\' el and J. Dhombres, Functional Equations in Several Vaiables, Cambridge Univ. Press, 1989.
\bibitem{CHO} P. W. Cholewa, Remarks on the stability of functional equations, Aequationes Math. 27 (1984), 76-86.
\bibitem{CZE1} S. Czerwik, On the stability of the quadratic mapping in normed spaces, Abh. Math. Sem. Univ. Hamburg 62 (1992), 59-64.
\bibitem{CZE2} S. Czerwik, Functional Equations and Inequalities in Several Variables, World Scientific, River Edge, NJ, 2002.
\bibitem{CZE3} S. Czerwik (ed.), Stability of Functional Equations of Ulam--Hyers--Rassias Type, Hadronic Press, 2003.
\bibitem{DRL} F. Drljevi\' c, On a functional which is quadratic on $A$-orthogonal vectors. Publ. Inst. Math. (Beograd)(N.S.) 54(1986), 63-71.
\bibitem{FOC} M. Fochi, Functional equations in $A$-orthogonal vectors, Aequationes Math. 38 (1989), 28-40.
\bibitem{G-S} R. Ger and J. Sikorska, Stability of the orthogonal additivity, Bull Polish Acad. Sci. Math. 43 (1995), No. 2, 143-151.
\bibitem{HYE} D. H. Hyers, On the stability of the linear functional equation, Proc. Nat. Acad.
Sci. U.S.A. 27 (1941), 222–224.
\bibitem{H-I-R} D. H. Hyers, G. Isac and Th. M. Rassias, Stability of Functional Equations in Several Variables, Birkh\"auser, Basel, 1998.
the Absence of the Mountain
\bibitem{J-L} K. W. Jun and Y. H. Lee, On the Hyers--Ulam--Rassias stability of a pexiderized
quadratic inequality, Math. Ineq. Appl., 4(1) (2001), 93–118.
\bibitem{J-S} S.-M. Jung and P. Sahoo, Hyers-Ulam stability of the quadratic equation of Pexider type, J. Korean Math. Soc. 38 (2001), No. 3, 645-656.
\bibitem{MOS} M. S. Moslehian, On the stability of the orthogonal Pexiderized Cauchy equation, arXiv math.FA/0412474.
\bibitem{RAS} Th.M. Rassias, On the stability of the quadratic functional equation and its applications, Studia Univ. Babe\c s-Bolyai Math. 43 (1998), no. 3, 89-124.
\bibitem{RAT} J. R\" atz, On orthogonally additive mappings, Aequations Math. 28 (1985), 35-49.
\bibitem{SKO} F. Skof , Local properties and approximations of operators, Rend. Sem. Mat. Fis. Milano 53 (1983) 113–129.
\bibitem{SZA} Gy. Szab\' o, Sesquilinear-orthogonally quadratic mappings, Aequationes Math. 40 (1990), 190-200.
\bibitem{VAJ} F. Vajzovi\' c, \" Uber das Funktional $H$ mit der Eigenschaft: $(x,y)=0 \Rightarrow H(x+y)+H(x-y)=2H(x)+2H(y)$, Glasnik Mat. Ser. III 2 (22)(1967), 73-81.
\end{thebibliography}
\end{document}